\documentclass[10pt]{article}

\usepackage{amsmath, amsthm, amssymb}
\usepackage{mathrsfs}
\usepackage{bm}
\usepackage{natbib}
\usepackage[dvips]{graphicx}

\newtheorem{thm}{Theorem}
\newtheorem{lem}{Lemma}

\newtheorem{definition}{Definition}
\newtheorem{remark}{Remark}
\newtheorem{example}{Example}

\newtheorem{exercise}{Exercise}
\newtheorem{report}{Report}
\newtheorem{assumption}{Assumption}

\newcommand{\Proof}{\noindent {\bf Proof: }\ }
\newcommand{\bin}{\mathrm{Bin}}

\newcommand{\prob}{\mathrm{Prob}}
\newcommand{\abs}[1]{\vert #1 \vert}
\newcommand{\abslr}[1]{\left\vert #1 \right\vert}

\newcommand{\real}{\mathbb{R}}
\newcommand{\dist}{\mathrm{dist}}

\newcommand{\rmd}{\mathrm{d}}
\newcommand{\rmor}{\;\mathrm{or}\;}
\newcommand{\io}{{i.o.}}
\newcommand{\locparam}{\mu}
\newcommand{\scparam}{\sigma}
\newcommand{\wtparam}{\alpha}

\newcommand{\scrg}{\mathscr{G}}
\newcommand{\scrk}{\mathscr{K}}

\newcommand{\seqnum}{j}

\title{Strong consistency of MLE for finite mixtures of location-scale distributions 
when the ratios of the scale parameters are exponentially small}

\author{Kentaro Tanaka\footnote{Department of Industrial Engineering and Management, Tokyo Institute of Technology, 
2-12-1 O-okayama, Meguro-ku, Tokyo 152-8552 JAPAN, E-mail:\,tanaka.k.al@m.titech.ac.jp}}

\begin{document}

\maketitle

\begin{abstract}
In finite mixtures of location-scale distributions, 
if there is no constraint on the parameters then the maximum likelihood estimate does not exist. 
But when the ratios of the scale parameters are restricted appropriately, 
the maximum likelihood estimate exists. 
We prove that the maximum likelihood estimator (MLE) is strongly consistent, 
if the ratios of the scale parameters are restricted from below by $\exp(-n^{d}) , \; 0 < d < 1 $, 
where $n$ is the sample size. 
\end{abstract}

\section{Introduction}
\label{sec:preliminaries}

In this paper we consider mixtures of $M$ location-scale densities 
which is defined by
\begin{eqnarray}
 f(x;\theta) \equiv \sum_{m=1}^{M} \wtparam_{m} f_{m}(x;\locparam_{m}, \scparam_{m}) , 
  \nonumber 
\end{eqnarray}
where 
$\wtparam_{m}$, called the mixing weights, 
are nonnegative real numbers that sum to one and  
$f_{m}(x;\locparam_{m}, \scparam_{m})$, called the components of the mixture,
are location-scale density functions with the location parameter $\locparam_{m}\in \mathbb{R}$
and the scale parameter $\scparam_{m}>0$. 
$\theta$ contains all the parameters in the mixture and can be written as 
$\theta = (\wtparam_{1}, \locparam_{1}, \scparam_{1}, \dots, \wtparam_{M}, \locparam_{M}, \scparam_{M})$. 

The location-scale densities satisfy 
\begin{equation}
f_{m}(x;\locparam_{m}, \scparam_{m})=\frac{1}{\scparam_m} f_m \left( \frac{x-\locparam_m}{\scparam_m}; 0,1\right).
\nonumber 
\end{equation}
We allow $f_{m}(x;\locparam_{m}, \scparam_{m})$ to belong to different families. 
For example, $f_{1}(x;\locparam_{1}, \scparam_{1})$ may be a normal density, 
$f_{2}(x;\locparam_{2}, \scparam_{2})$ may be a uniform density, etc.

In finite mixtures of location-scale distributions, 
the maximum likelihood estimate over the whole parameter space does not exist. 
For a given data $x_{1},\dots,x_{n}$, the likelihood function for a mixture is unbounded. 
For example, if we set $\locparam_{1}=x_{1}$ and let $\scparam_{1}\rightarrow 0$ in a mixture of two normal densities 
with means $\locparam_{1}, \locparam_{2}$ and standard deviations $\scparam_{1}, \scparam_{2}$, 
then the likelihood will tend to infinity. 

Let us consider the constrained parameter space
\begin{equation}
 \Theta_{b} \equiv \{\theta \in \Theta \mid \min_{1 \leq m,m' \leq M}\scparam_{m}/\scparam_{m'} \geq b\}
  \qquad , \qquad 0 < b < 1, 
 \nonumber 
\end{equation}
where $\Theta$ denote the unconstrained parameter space. 
\citet{H1985} showed that the global maximizer $\hat{\theta}_{b}$ of the likelihood function 
over $\Theta_{b}$ exists and if true value of parameter belongs to $\Theta_{b}$ 
then $\hat{\theta}_{b}$ is strongly consistent. 
But there is the problem how small we have to choose $b$ to ensure strong consistency. 
An interesting question here is whether we can decrease the bound $b$ to zero with the 
sample size yet gurantee the strong consistency of the maximum likelihood estimator. 
This question is mentiond in \citet{H1985}, \citet{M2000} and treated as an unsolved problem. 

Meanwhile in \citet{TT2003-35}, we consider mixtures of location-scale distributions with constraint imposed on 
the scale parameters themselves and showed that 
the maximum likelihood estimator is strongly consistent 
if the scale parameters are restricted from below by $\exp(n^{-d}),\; 0<d<1$. 
\citet{TT2003-20} implies that the rate $\exp(n^{-d})$ obtained in \citet{TT2003-35} is almost 
the lower bound to ensure strong consistency. 
The method used in \citet{TT2003-35} is useful for solving the problem stated in \citet{H1985} 
in which the constraints are imposed on the ratios of the scale parameters. 

In this paper, we solve the problem stated in \citet{H1985}. 
We prove that the maximum likelihood estimator is strongly consistent, 
if the ratios of the scale parameters are restricted from below by $\exp(-n^{d}) , \; 0 < d < 1 $. 

The organization of the paper is as follows. 
In section~\ref{sec:notation-definitions} we prepare notation and summarize some preliminary results. 
In section~\ref{sec:main-result} we state our main result in theorem~\ref{thm:bn-thm}. 
Section~\ref{sec:proof} is devoted to the proof of theorem~\ref{thm:bn-thm}. 
The last section is conclusion and future work. 

\section{Notation and definitions}
\label{sec:notation-definitions}
\subsection{Notation}

Let 
$\Omega_{m}=\mathbb{R}\times (0,\infty)$
denote the parameter space of the $m$-th component $(\locparam_m, \scparam_m)$. 
Then the entire parameter space $\Theta$ can be represented as follows.
\begin{eqnarray}
 \Theta  & = & 
 \{
 (\wtparam_{1}, \ldots , \wtparam_{M}) \in \real^{M}
 \mid 
 \sum_{m=1}^{M} \wtparam_{m} = 1 
 \; , \;
 \wtparam_{m} \geq 0
 \}\times \prod_{m=1}^{M}\Omega_{m}.
 \nonumber 
\end{eqnarray}

Let $\scrk$ be a subset of $\{1,2,\ldots,M\}$ and 
let $\abslr{\scrk}$ denote the number of elements in
$\scrk$. 
Denote by $\theta_{\scrk}$ a subvector of $\theta \in \Theta$ 
consisting of the components in $\scrk$. 
Then the parameter space of subprobability measures 
consisting of the components in $\scrk$ 
is
\begin{equation}
 \bar{\Theta}_{\scrk} \equiv 
  \{
  \theta_{\scrk} \mid \theta \in \Theta, \sum_{m \in \scrk} \wtparam_m
  \le 1 
  \}.
  \nonumber 
\end{equation}
Corresponding density and the set of subprobability densities are
denoted by 
\begin{align}
 f_{\scrk}(x;\theta_{\scrk})
 & \equiv  \sum_{m \in \scrk} \wtparam_{m} f_m(x;{\locparam}_m, \scparam_{m}) ,
 \nonumber 
 \\
 \scrg_{\scrk} 
 & \equiv 
 \{
 f_{\scrk}(x;\theta_{\scrk}) 
 \mid 
 \theta_{\scrk} \in \bar{\Theta}_{\scrk}
 \} .
 \nonumber 
\end{align}
Furthermore denote 
the set of subprobability densities with no more than $K$
components by
\begin{equation}
 \scrg_{K} \equiv \bigcup_{\abslr{\scrk}\leq K}\scrg_{\scrk}
  \qquad (1 \leq K \leq M). 
  \nonumber 
\end{equation}

\subsection{Identifiability and strong consistency of estimators}

In general, a parametric family of distributions is identifiable if
different values of parameter designate different distributions.  In
mixtures of distributions, different parameters may designate the same
distribution.  For example, if $\wtparam_{1} = 0$, then for all
parameters which differ only in $\locparam_{1}$ or $\scparam_{1}$, we have the same distribution. 
Thus mixtures of distributions are not identifiable.
Therefore we have to carefully define strong consistency of an estimator.  
The following definition is essentially the same as \citet{R1981}.  
We assume that the parameter space $\Theta$ is a subset of Euclidean space and
$\dist(\theta,\theta')$ denotes the Euclidean distance between
$\theta,\theta' \in \Theta$.
Furthermore we define  
\begin{eqnarray}
 \dist(U,V) 
  \equiv 
  \inf_{\theta \in U} \inf_{\theta' \in V} \dist(\theta,\theta')
  \nonumber 
\end{eqnarray}
for $U,V \subset \Theta$.

\begin{definition} 
\label{definition:1}
 $($a strongly consistent estimator$)$ \\
 Let $\Theta_0$ denote the set of true parameters
 \begin{eqnarray}
  \Theta_0 \equiv \{\theta \in \Theta \mid f(x;\theta) = f(x;\theta_0)
   \quad a.e. \; x\} , 
   \nonumber
 \end{eqnarray}
 where $\theta_{0}$ is one of parameters 
 designating the true  distribution, 
 and let 
 \begin{eqnarray}
  \Theta(\hat{\theta}) 
   \equiv 
   \{\theta \in \Theta \mid f(x;\theta) = f(x;\hat{\theta})
   \quad a.e. \; x\} . 
   \nonumber
 \end{eqnarray}
 An estimator $\hat{\theta}_n$ is strongly consistent if
  \begin{eqnarray}
  \prob
   \left(
    \lim_{n \rightarrow \infty}
    \dist(\Theta(\hat{\theta}_{n}), \Theta_0) = 0
   \right) = 1 . 
   \nonumber
 \end{eqnarray}
\end{definition}

\subsection{Preliminaries}

In \citet{TT2003-35}, 
we assume the following regularity conditions for strong consistency of the constrained maximum likelihoood estimator.  
\begin{assumption}
 There exist real constants 
 $v_{0},v_{1} > 0$ and $\beta > 1$ such that 
 \begin{eqnarray}
  f_{m}(x;\locparam_{m}=0, \scparam_{m}=1) 
   \leq 
   \min \{v_{0} \;, \; v_{1} \cdot \abs{x}^{-\beta}\}
   \nonumber 
 \end{eqnarray}
 for all $m$.
\label{assumption:1}
\end{assumption}
Let $\Gamma$ denote any compact subset of $\Theta$. 
\begin{assumption}
 For $\theta \in \Theta$ and any positive real number $\rho$, let 
 \begin{eqnarray}
  f(x;\theta,\rho) 
   & \equiv & \sup_{\dist(\theta',\theta) \leq \rho}f(x;\theta')
   . 
   \nonumber
 \end{eqnarray}
 For each $\theta \in \Gamma$ and sufficiently small $\rho$, 
 $f(x;\theta, \rho)$ is measurable. 
\label{assumption:2}
\end{assumption}
\begin{assumption}
 For each $\theta \in \Gamma$,  if\/
 $\lim_{\seqnum \rightarrow \infty}\theta^{(\seqnum)} = \theta, \; (\theta^{(\seqnum)} \in \Gamma)$ 
 then\/
 $\lim_{\seqnum \rightarrow \infty}
 f(x;\theta^{(\seqnum)}) = f(x;\theta)$
 except on a set which is a null set and 
 does not depend on the sequence 
 $\{\theta^{(\seqnum)}\}_{\seqnum=1}^{\infty}$.
\label{assumption:3}
\end{assumption}
\begin{assumption}
 \begin{eqnarray}
  \int \abslr{
   \log{f(x;\theta_0)}
   }
   f(x;\theta_0)\rmd x < \infty .
   \nonumber 
 \end{eqnarray}
\label{assumption:4}
\end{assumption}

Let $E_{0}[\cdot]$ denote the expectation under the true parameter $\theta_{0}$. 
In \citet{TT2003-35}, we showed the following theorem. 
\begin{thm}{\rm (\citet{TT2003-35})}
 Suppose that assumptions \ref{assumption:1}--\ref{assumption:4}
 are satisfied 
 and $f(x;\theta_{0}) \in \scrg_{M}\backslash \scrg_{M-1}$. 
 Let $c_{0} > 0$ and $ 0 < d < 1$.
 If $c_{n} = c_{0} \cdot \exp(-n^{d})$ and 
 \begin{eqnarray}
  \Theta_{c_{n}} 
   \equiv 
   \{
   \theta \in \Theta
   \mid 
   \min_{1\leq m \leq M} \scparam_{m} \geq c_{n}
   \},
   \nonumber 
 \end{eqnarray}
 then
  \begin{eqnarray}
  \prob
   \left(
    \lim_{n \rightarrow \infty}
    \dist(\Theta(\hat{\theta}_{c_{n}}), \Theta_0) = 0
   \right) = 1 \; ,  
   \nonumber 
\end{eqnarray}
where $\hat{\theta}_{c_{n}}$ is the maximum likelihood estimator 
restricted to $\Theta_{c_{n}}$.  
 \label{thm:cn-thm}
\end{thm}

\section{Main result}
\label{sec:main-result}

To show the strong consistency of the constrained maximum likelihoood estimator
in the problem stated in \citet{H1985}, 
we replace the assumption~\ref{assumption:1} with the following assumption. 
\begin{assumption}
 There exist real constants 
 $v_{0},v_{1} > 0$ and $\beta > 2$ such that 
 \begin{eqnarray}
  f_{m}(x;\locparam_{m}=0, \scparam_{m}=1) 
   \leq 
   \min \{v_{0} \;, \; v_{1} \cdot \abs{x}^{-\beta}\}
   \nonumber 
 \end{eqnarray}
 for all $m$.
\label{assumption:5}
\end{assumption}

Now we state the main theorem of this paper. 
\begin{thm}
 \label{thm:bn-thm}
 Suppose that the assumptions~\ref{assumption:2}--\ref{assumption:5} are satisfied 
 and $f(x;\theta_{0}) \in \scrg_{M}\backslash \scrg_{M-1}$. 
 Let $b_{0} > 0$ and $0 < d < 1$. 
 If $b_{n} = b_{0} \cdot \exp(-n^{d})$ 
 and 
  \begin{eqnarray}
  \Theta_{b_{n}} 
   \equiv 
   \{
   \theta \in \Theta
   \mid 
   \min_{1\leq m\neq m' \leq M}\frac{\scparam_{m}}{\scparam_{m'}} \geq b_{n}
   \}, 
   \nonumber 
 \end{eqnarray}
 then 
  \begin{eqnarray}
  \prob
   \left(
    \lim_{n \rightarrow \infty}
    \dist(\Theta(\hat{\theta}_{b_{n}}), \Theta_0) = 0
   \right) = 1 \; ,  
   \nonumber 
\end{eqnarray}
where $\hat{\theta}_{b_{n}}$ is the maximum likelihood estimator restricted to $\Theta_{b_{n}}$. 
\qed
\end{thm}

\section{Proof}
\label{sec:proof}

In this section, we prove theorem~\ref{thm:bn-thm} by using theorem~\ref{thm:cn-thm}. 
The organization of this section is as follows. 
First in subsection~\ref{sec:part-param-space} we partition the parameter space $\Theta_{b_{n}}$ into two sets. 
Then the proof for strong consistency of the maximum likelihood estimator restricted to $\Theta_{b_{n}}$
is also partitioned and the proof for one set is shown by applying the result of theorem~\ref{thm:cn-thm}. 
The proof for another set is shown in section~\ref{sec:proof-lemma-refl}. 

\subsection{Partitioning the parameter space}
\label{sec:part-param-space}

Let $0 < d < 1$ be any constant and define $b_{n} \equiv \exp(-n^{d})$. 
We choose $d'$ such that $d < d' < 1$ and define $c_{n} \equiv \exp(-n^{d'})$. 
Notice that the following arguments also hold
even when we define $b_{n} \equiv b_{0}\cdot\exp(-n^{d})$ with a positive constant $b_{0}$. 
Define 
$\Theta_{b_{n}}\equiv\{\theta \in \Theta \mid \min_{1\leq m\neq m' \leq M}\frac{\scparam_{m}}{\scparam_{m'}} \geq b_{n}\}$
and 
$\Theta_{c_{n}}\equiv\{\theta \in \Theta \mid \min_{1\leq m \leq M}{\scparam_{m}} \geq c_{n}\}$. 
The constrained parmeter space $\Theta_{b_{n}}$ can be divided into two sets.  
\begin{equation}
 \Theta_{b_{n}} = (\Theta_{b_{n}} \cap \Theta_{c_{n}}) \cup (\Theta_{b_{n}} \cap \Theta_{c_{n}}^{C}), 
  \nonumber 
\end{equation}
where $\Theta_{c_{n}}^{C}$ is the complement of $\Theta_{c_{n}}$. 

From theorem~\ref{thm:cn-thm}, the maximum likelihood estimator over $\Theta_{b_{n}} \cap \Theta_{c_{n}}$ 
is strongly consistent. 
If the maximum of likelihood function over $\Theta_{b_{n}} \cap \Theta_{c_{n}}^{C}$ is very small, 
then the maximum likelihood estimator over $\Theta_{b_{n}}$ is strongly consisitent. 
By 
the argument used in \citet{W1949}, in order to prove theorem~\ref{thm:bn-thm} 
it suffices to prove the following lemma. 
\begin{lem}
 \label{lem:goal}
 \begin{equation}
  \lim_{n \rightarrow \infty}
  \frac{
   \sup_{\theta \in \Theta_{b_n}\cap\Theta_{c_n}^{C}}
   \prod_{i=1}^{n} f(x_{i};\theta)
   }
   { \prod_{i=1}^{n} f(x_{i};\theta_{0}) } = 0, \quad a.e. 
   \label{eq:lem:goal}
 \end{equation}
\end{lem}

\subsection{Proof of lemma~\ref{lem:goal}}
\label{sec:proof-lemma-refl}

Let 
\begin{equation}
 \scparam_{(1)} = \min_{1\leq m \leq M} \sigma_{m}
  \qquad , \qquad
 \scparam_{(M)} = \max_{1\leq m \leq M} \sigma_{m}. 
 \label{eq:OrderedScparam}
\end{equation}
Then for $\theta \in \Theta_{c_n}^{C}$, 
\begin{equation}
 \scparam_{(1)} \leq c_{n}. 
  \label{eq:4} 
\end{equation}
Furthermore for $\theta \in \Theta_{b_n}$, 
\begin{equation}
 \frac{\scparam_{(1)}}{\scparam_m} \geq b_{n}
  \quad , \quad 
  1 \leq m \leq M. 
  \nonumber 
\end{equation}
Therefore for $\theta \in \Theta_{b_n}\cap\Theta_{c_n}^{C}$, 
\begin{equation}
 \scparam_m \leq c_{n}/b_{n} = \exp(n^{d}-n^{d'})  
  \quad , \quad 
  1 \leq m \leq M. 
  \label{eq:BoundingTheScaleParameter}
\end{equation}
This means that all the scale paramters $\scparam_{m}$ of $\theta \in \Theta_{b_n}\cap\Theta_{c_n}^{C}$ are very small for large $n$. 
Hence the maximum of likelihood function over $\Theta_{b_n}\cap\Theta_{c_n}^{C}$ seems to be small. 

Let $E_{0}[\cdot]$ denote the expectation under the true parameter $\theta_{0}$.
By law of large numbers (\ref{eq:lem:goal}) is implied by 
\begin{equation}
 \limsup_{n\rightarrow \infty} \frac{1}{n} \sup_{\theta \in \Theta_{b_n}\cap\Theta_{c_n}^{C}}
  \sum_{i=1}^{n}\log{f(x_{i};\theta)} < E_{0}[\log{f(x;\theta_{0})}] \quad a.e. 
  \label{eq:goal2}
\end{equation}
Therefore, in order to prove~\ref{lem:goal} it suffices to prove (\ref{eq:goal2}). 

\subsubsection{Step 1 : Bounding the components by step functions}
\begin{figure}[htbp]
 \begin{center}
  \includegraphics[width=8cm]{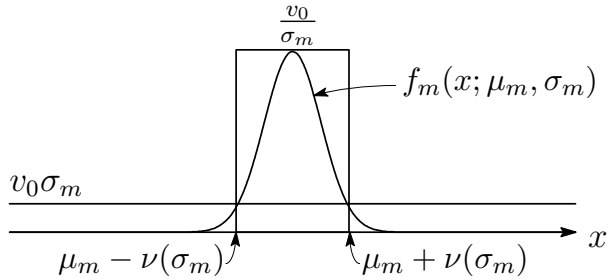}
  \caption{Each component is bounded by a step function.}
  \label{fig:BoundedByStepFunction}
 \end{center}
\end{figure}
Define 
\begin{equation}
 \nu(\scparam) \equiv \left(\frac{v_{1}}{v_{0}}\right)^{\frac{1}{\beta}}\cdot\scparam^{1-\frac{2}{\beta}}. 
  \label{eq:DefinitionOfNu} 
\end{equation}
From assumption~\ref{assumption:5} each component is bounded from above as 
\begin{equation}
 f_{m}(x;\locparam_{m}, \scparam_{m}) \leq \max\{1_{[\locparam_{m}-\nu(\scparam_{m}),\; \locparam_{m}+\nu(\scparam_{m})]}(x)\cdot\frac{v_{0}}{\scparam_{m}}
  \; , \; {v_{0}}{\scparam_{m}}\}. 
  \nonumber 
\end{equation}
See figure~\ref{fig:BoundedByStepFunction}. 
From this and (\ref{eq:OrderedScparam}) we obtain the following lemma. 
\begin{lem}
 \label{lem:BoundingTheComponents}
 \begin{equation}
   f_{m}(x;\locparam_{m}, \scparam_{m}) \leq \max\{1_{[\locparam_{m}-\nu(\scparam_{m}),\; \locparam_{m}+\nu(\scparam_{m})]}(x)\cdot\frac{v_{0}}{\scparam_{(1)}}
   \; , \; {v_{0}}{\scparam_{(M)}}\}
   \quad , \quad 
   1 \leq m \leq M
   . 
   \nonumber 
 \end{equation}
\end{lem}

\subsubsection{Step 2 : Bounding the likelihood function by two terms}

Let $R_{n}(V)$ denote the number of observation which belong to a set $V \subset \real$. 
Define 
\begin{equation}
 J(\theta) \equiv \bigcup_{m=1}^{M}[\locparam_{m}-\nu(\scparam_{m}), \; \locparam_{m}+\nu(\scparam_{m})]. 
  \nonumber 
\end{equation}
\begin{lem}
\begin{equation}
 \forall \theta \in \Theta_{b_{n}}\cap\Theta_{c_{n}}^{C}
  \quad , \quad 
 \sum_{i=1}^{n}\log{f(x_{i};\theta)} \leq R_{n}(J(\theta))\cdot\log{\frac{v_{0}}{\scparam_{(1)}}}
 + R_{n}(J(\theta)^{C})\cdot\log{(v_{0}\scparam_{(M)})}
  \label{eq:BoundingTheLikelihoodFunction}
\end{equation}
\qed
\end{lem}
\Proof
From lemma~\ref{lem:BoundingTheComponents}, we obtain 
\begin{eqnarray}
 \sum_{i=1}^{n}\log{f(x_{i};\theta)}
  & = &
  \sum_{i=1}^{n}\log{\left\{\sum_{m=1}^{M}\wtparam_{m}f_{m}(x_{i};\locparam_{m}, \scparam_{m})\right\}}
 \nonumber \\ & \leq &
 \sum_{i=1}^{n}\left\{\max_{m=1,\dots,M}\log{f_{m}(x_{i};\locparam_{m}, \scparam_{m})}\right\} 
 \nonumber \\ & \leq &
 \sum_{i=1}^{n}
 \max_{m=1,\dots,M} \max\{
 1_{[\locparam_{m}-\nu(\scparam_{m}), \; \locparam_{m}+\nu(\scparam_{m})]}(x) \cdot \log{\frac{v_{0}}{\scparam_{(1)}}}  
 \; , \; 
 \log{(v_{0}\scparam_{(M)})}
 \}
 \nonumber \\ & = & 
 R_{n}(J(\theta))\cdot\log{\frac{v_{0}}{\scparam_{(1)}}} + R_{n}(J(\theta)^{C})\cdot\log{(v_{0}\scparam_{(M)})}. 
 \nonumber 
\end{eqnarray}
\qed

In the following we bound the right hand side of (\ref{eq:BoundingTheLikelihoodFunction}) from above. 

\subsubsection{Step 3 : Bounding the first term}

Let  $x_{1}, \ldots, x_{n}$ denote a random sample 
of size $n$ from $f(x;\theta_{0})$ and let 
\begin{eqnarray}
 x_{n,1}  \equiv  \min{\{x_{1}, \ldots, x_{n}\}}
  \quad , \quad 
 x_{n,n}  \equiv  \max{\{x_{1}, \ldots, x_{n}\}}.
  \nonumber 
\end{eqnarray}
In \citet{TT2003-35}, we showed the following lemma. 
\begin{lem}
 \label{lem:BoundedInterval}
 {\rm(\citet{TT2003-35})}
 For any real positive constants $A_{0} > 0 , \zeta > 0$, define 
 \begin{equation}
  A_{n} = A_{0}\cdot n^{\frac{2+\zeta}{\beta-1}}. 
   \label{eq:An} 
 \end{equation}
 Then
 \begin{equation}
  \prob\left(
	x_{n,1} < -A_{n} \; \rmor \;  x_{n,n} > A_{n}
	\quad \io
       \right) = 0. 
  \nonumber 
 \end{equation}
\end{lem}
By this lemma we can bound
the behavior of the minimum and the maximum of the sample with probability 1. 
In the following we ignore the event $\{x_{n,1} < -A_{n} \; \rmor \;  x_{n,n} > A_{n}\}$. 

Next we prove the following lemma for bounding the first term of (\ref{eq:BoundingTheLikelihoodFunction}). 
\begin{lem}
 \label{lem:TheNumberOfTheObservationsWithinShortInterval}
 \begin{equation}
  \forall \theta \in \Theta_{b_{n}}\cap\Theta_{c_{n}}^{C}
   \; , \; 
  \forall \epsilon > 0
     \quad , \quad 
  \prob\left(
	\max\{R_{n}(J(\theta))-4M \; ,\; 0\} > \epsilon \quad \io
       \right) = 0 
  \nonumber
 \end{equation}
\end{lem}
\Proof
\begin{figure}[htbp]
 \begin{center}
  \includegraphics[width=10cm]{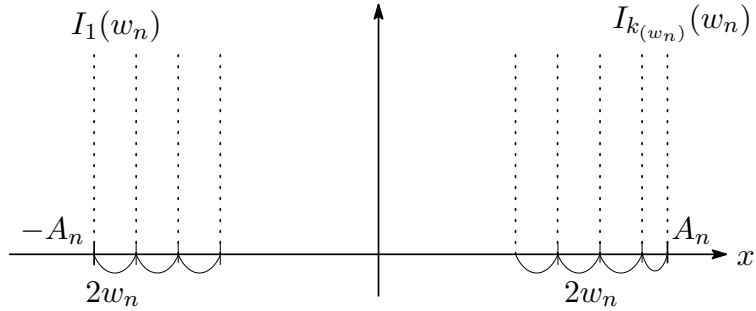}
  \caption{Division of interval $[-A_{n}, A_{n}]$ by short interval of length $2w_{n}$.}
  \label{fig:ShortIntervals}
 \end{center}
\end{figure}
Define 
\begin{equation}
 w_{n}\equiv \nu(c_{n}/b_{n}). 
 \label{eq:wn}
\end{equation}
From lemma~\ref{lem:BoundedInterval}, we can ignore the event $\{x_{n,1} < -A_{n} \; \rmor \;  x_{n,n} > A_{n}\}$. 
Now we divide $[-A_{n},A_{n}]$ from $-A_{n}$ to $A_{n}$ by short intervals of length $2w_{n}$. 
Let $k(w_{n})$  be the number of the short intervals and let 
$I_{1}(w_{n}),\dots,I_{k(w_{n})}(w_{n})$
be the divided short intervals. 
The length of the rightmost short interval $I_{k(w_{n})}(w_{n})$ may be less than $2w_{n}$. 
See figure~\ref{fig:ShortIntervals}. 
Then we have 
\begin{equation}
 k(w_{n}) \leq \frac{2A_{n}}{2w_{n}} + 1 = \frac{A_{n}}{w_{n}} + 1. 
  \label{eq:6} 
\end{equation} 
From (\ref{eq:BoundingTheScaleParameter}), (\ref{eq:DefinitionOfNu}) and (\ref{eq:wn}) we have
\begin{equation}
 \nu(\scparam_1) ,\; \nu(\scparam_2) ,\; \dots ,\; \nu(\scparam_M) \leq \nu(c_{n}/b_{n}) = w_{n}. 
  \nonumber 
\end{equation}
Since $J(\theta)=\bigcup_{m=1}^{M}[\locparam_{m}-\nu(\scparam_{m}), \; \locparam_{m}+\nu(\scparam_{m})]$ 
consists of M intervals of length at most $2w_{n}$, 
$J(\theta)\cap [-A_{n},A_{n}]$ is 
covered by at most $2M$ short intervals of $I_{1}(w_{n}),\dots,I_{k(w_{n})}(w_{n})$. 
Therefore the following relation holds. 
\begin{eqnarray}
 \{\max\{R_{n}(J(\theta))-4M \; ,\; 0\} > \epsilon \}
  & \Leftrightarrow & 
 \{\max\{R_{n}(J(\theta)\cap [-A_{n},A_{n}])-4M \; ,\; 0\} > \epsilon \}
 \nonumber \\ 
  & \Rightarrow &
  \{1\leq\exists k \leq k(w_{n}) \; , \; R_{n}(I_{k}(w_{n})) \geq 2\}.  
  \nonumber 
\end{eqnarray}
Let $u_{0} \equiv \sup_{x}f(x;\theta_{0})$. 
Let $P_0(V)$ denote the probability of
$V \subset \real$ under the true density
\begin{eqnarray}
 P_0(V) \equiv \int_{V} f(x;\theta_{0}) \rmd x \; . 
  \nonumber 
\end{eqnarray}
From (\ref{eq:6}), $R_n(I_{k}(w_{n})) \sim \bin(n , P_0(I_{k}(w_{n})))$ and $P_0(I_{k}(w_{n})) \le 2w_{n} u_{0}$
we obtain 
\begin{eqnarray}
 \lefteqn{
  \prob
  \left(
   \max\{R_{n}(J(\theta))-4M \; ,\; 0\} > \epsilon
  \right)
  \leq
  \sum_{k=1}^{k(w_{n})}
  \prob
  \left(
   R_n(I_{k}(w_{n})) \geq 2
  \right)
  } \hspace{2cm} & & 
 \nonumber \\ & \leq &
  k(w_{n})
  \cdot 
  \left\{
  \max_{1\le k \le k(w_{n})}
  \prob(R_n(I_{k}(w_{n})) \geq 2)
  \right\}
 \nonumber \\ & \leq & 
  \left(
   \frac{A_{n}}{w_{n}} + 1
  \right)
  \sum_{k=2}^{n}
  \binom{n}{k}
  (2w_{n} u_{0})^{k}(1 - 2w_{n} u_{0})^{n-k}
  \nonumber \\ & \leq &
  \left(
   \frac{A_{n}}{w_{n}} + 1
  \right)
  \sum_{k=2}^{n}\frac{n^{k}}{k!}(2w_{n} u_{0})^{k}
  \nonumber \\ & \leq &
    \left(
   \frac{A_{n}}{w_{n}} + 1
  \right)
  (2nw_{n} u_{0})^{2}
  \exp{(2nw_{n} u_{0})}
  \; .
  \label{eq:3}
\end{eqnarray}
From 
(\ref{eq:DefinitionOfNu}), (\ref{eq:An}) and (\ref{eq:wn}), 
the order of the right hand side of (\ref{eq:3}) is evaluated as follws. 
\begin{eqnarray}
 \left(\frac{A_{n}}{w_{n}} + 1\right)(2nw_{n} u_{0})^{2}
  & = &
  O\left(
    n^{2+\frac{2+\zeta}{\beta-1}}\cdot e^{(n^{d}-n^{d'})(1-2/\beta)}
   \right)
  \nonumber \\ 
 \exp{(2nw_{n} u_{0})} & = & O(1)
  \nonumber 
\end{eqnarray}
Recall that $0 < d < d' < 1$ and $\beta > 2$ by assumption~\ref{assumption:5}. 
Then 
\begin{equation}
 \sum_{n=1}^{\infty} n^{2+\frac{2+\zeta}{\beta-1}}\cdot e^{(n^{d}-n^{d'})(1-2/\beta)} < \infty. 
  \nonumber 
\end{equation}
Therefore, when we sum the right hand side of (\ref{eq:3}) over $n$, the resulting series converges. 
Hence by Borel-Cantelli lemma, we have 
 \begin{equation}
  \prob\left(
	\max\{R_{n}(J(\theta))-4M \; ,\; 0\} > \epsilon \quad \io
       \right) = 0. 
  \nonumber
 \end{equation}
\qed

\subsubsection{Step 4 : Evaluating the likelihood function}

From lemma~\ref{lem:TheNumberOfTheObservationsWithinShortInterval}, 
we  can ignore the event $R_{n}(J(\theta)) > 4M$. 
In the following we consider only the event $\{R_{n}(J(\theta)) \leq 4M\}$. 
Note that $R_{n}(J(\theta)^{C}) \geq n-4M$. 
For sufficiently large $n$, 
the scale parameters $\sigma_{1},\dots,\sigma_{M}$ of $\theta \in \Theta_{b_{n}}\cap\Theta_{c_{n}}^{C}$ are very small so that 
we can bound the right hand side of (\ref{eq:BoundingTheLikelihoodFunction}) as follows. 
\begin{eqnarray}
 \lefteqn{
 R_{n}(J(\theta))\cdot\log{\frac{v_{0}}{\scparam_{(1)}}} + R_{n}(J(\theta)^{C})\cdot\log{(v_{0}\scparam_{(M)})}
  } \hspace{2cm} & & 
  \nonumber \\ 
  & \leq &
  4M\cdot\log{\frac{v_{0}}{\scparam_{(1)}}} + (n-4M)\cdot\log{(v_{0}\scparam_{(M)})} 
  \qquad 
  \nonumber \\ & \leq &
  n\cdot\log{v_{0}} + (n-8M)\cdot \log{\scparam_{(1)}} + (n-4M)\cdot\log{\frac{1}{b_{n}}} 
    \quad a.e. 
    \label{eq:5}
\end{eqnarray}
where the last inequality holds by (\ref{eq:BoundingTheScaleParameter}) i.e. $\scparam_{(M)}\leq\scparam_{(1)}/b_{n}$. 
Recall that we set $0 < d < d' < 1$. 
Then from (\ref{eq:4}), (\ref{eq:BoundingTheLikelihoodFunction}) and (\ref{eq:5}) we have 
\begin{eqnarray}
 & & \hspace{-1.7cm} \sup_{\theta\in \Theta_{b_{n}}\cap \Theta_{c_{n}}^{C}}\frac{1}{n}\sum_{i=1}^{n}\log{f(x_{i};\theta)}
 \nonumber \\ & \leq &
 \frac{1}{n}
 \left\{
  n\cdot\log{v_{0}} + (n-8M)\cdot (-n^{d'}) + (n-4M)\cdot n^{d}
 \right\} 
 \rightarrow
  -\infty \quad a.e. 
 \nonumber 
\end{eqnarray}
Therefore we obtain (\ref{eq:goal2}) and lemma~\ref{lem:goal} is proved. 

This completes the proof of theorem~\ref{thm:bn-thm}. 

\section{Conclusion}
\label{sec:conclusion}

In this paper we prove that if we set $b_{n}\equiv \exp(-n^{d}), \;0 < d < 1$, then 
the maximum likelihood estimator restricted to 
$\Theta_{b_{n}}\equiv\{\theta \in \Theta \mid \min_{1\leq m\neq m' \leq M}\frac{\scparam_{m}}{\scparam_{m'}} \geq b_{n}\}$ 
is strongly consistent under very mild regularity conditions. 
Mixtures of normal distributions satisfy the regularity conditions. 
This means that the problem stated in \citet{H1985} is solved. 

If we define $b_{n}\equiv \exp(-n^{r}), \;r > 1$, and set $\theta$ as 
\begin{eqnarray}
 \locparam_{1}=x_{1} \quad &,& \quad \scparam_{1}=\exp(-n^{r})
  \nonumber \\ 
  \locparam_{m}=0 \quad &,& \quad \scparam_{m}=1 \qquad\qquad\quad (m\neq 1), 
   \nonumber 
\end{eqnarray}
then $\theta \in \Theta_{b_{n}}$ and the mean log likelihood of this density tends to infinity (\citet{TT2003-20}). 
This means that the mean log likelihood of the true model which converges to finite value almost everywhere is 
dominated by that of other models. 
Therefore if $b_{n}$ decreases to zero faster than $\exp(-n)$, 
then the consistency of the maximum likelihood estimator fails. 
This implies that the rate of $b_{n}\equiv \exp(-n^{d}), \;0 < d < 1$ 
obtained in this paper is almost the lower bound of the order of $b_{n}$ which maintains the consistency. 

In theorem~\ref{thm:bn-thm} we assume 
$f(x;\theta_{0}) \in \scrg_{M}\backslash \scrg_{M-1}$ 
i.e. the number of components of true model is known. 
To discuss the case that the number of components of true model is unknown, 
more complicated mathematical techniques are needed.

\bibliography{mixture}

\end{document}